\documentclass[a4paper,10pt,twoside,reqno]{amsart}
\usepackage{amsmath,amsfonts,amssymb,amsxtra,setspace,xspace,graphicx,lmodern,psfrag,epsfig,color,latexsym}
\usepackage[colorlinks=true]{hyperref}\hypersetup{urlcolor=blue, citecolor=red}

\newtheorem{thm}{Theorem}[section]
\newtheorem{cj}{Conjecture}[section]
\newtheorem{defi}[thm]{Definition}
\newtheorem{lem}[thm]{Lemma}
\newtheorem{cor}[thm]{Corollary}
\newtheorem{rem}[thm]{Remark}
\newtheorem{prop}[thm]{Proposition}

\newtheorem{assu}[thm]{Assumption}

\def\bthm{\begin{thm}\def\ethm{\end{thm}}}

\def\bdefi{\begin{defi}\sl{}\def\edefi{\end{defi}}}
\def\blem{\begin{lem}\def\elem{\end{lem}}}
\def\bcor{\begin{cor}\def\ecor{\end{cor}}}
\def\bprop{\begin{prop}\def\eprop{\end{prop}}}
\def\bpf{\begin{proof}}\def\epf{\end{proof}}

\newcommand{\R}{{\mathbb R}}

\newcommand{\be}[1]{\begin{equation}\label{#1}}
\newcommand{\ee}{\end{equation}}

\def\bb{\mathbb}
\def\ca{\mathcal}

\def\e{\varepsilon}

\def\bbd{\mathbf}
\def\E{\bb{E}}
\def\indiq{\bbd{1}}

\newcommand{\nc}{\normalcolor}

\begin{document}

\title[Jumping SDE\MakeLowercase{s} and nonlocal PDE\MakeLowercase{s}]
{On the equivalence between some jumping SDE\MakeLowercase{s}
with rough coefficients and some non-local PDE\MakeLowercase{s} }
\author{Nicolas Fournier}
\author{Liping Xu}

\address{N. Fournier, L. Xu: Laboratoire de Probabilit\'es et Mod\`eles al\'eatoires,
UMR 7599, UPMC, Case 188,
4 pl. Jussieu, F-75252 Paris Cedex 5, France.}
\email{nicolas.fournier@upmc.fr, liping.xu@upmc.fr}



\thispagestyle{empty}
\parindent=0pt
\parskip=6.3pt

\begin{abstract}
We study some jumping SDE
and the corresponding Fokker-Planck (or Kolmogorov forward) equation, which is a non-local PDE.
We assume only some measurability and growth conditions on the coefficients.
We prove that for any weak solution $(f_t)_{t\in [0,T]}$ of the PDE,
there exists a weak solution to the SDE of which the time marginals are given by $(f_t)_{t\in[0,T]}$.
As a corollary, we deduce that for any given initial condition,
existence for the PDE is equivalent to weak existence for the SDE
and uniqueness in law for the SDE implies uniqueness for the PDE.
This extends some ideas of Figalli \cite{MR2375067} concerning continuous SDEs and local PDEs.
\end{abstract}
\maketitle


\section{Introduction}

We consider the $d$-dimensional stochastic differential equation posed on some time interval $[0,T]$
\be{SDE}
X_t=X_0+\int_0^t b(s,X_s)\, ds+\int_0^t \sigma(s,X_s)\, dB_s+\int_0^t\int_{E} h(s,z,X_{s-})\, N(ds,dz),
\ee
where $(B_t)_{t\in[0,T]}$ is a $d$-dimensional Brownian motion and $N(ds,dz)$ is a Poisson measure
on $[0,T]\times E$
with intensity measure $ds\,\mu(dz)$.
The coefficients $b:[0,T]\times \R^d\mapsto \R^d$, $\sigma: [0,T]\times \R^d\mapsto\ca{S}_d^+$
and $h:[0,T]\times E \times \R^d  \mapsto \R^d$ are supposed to be at least measurable.
The space $E$ is endowed with a $\sigma$-field $\ca E$ and with a $\sigma$-finite measure
$\mu$ and $\ca S_d^+$ is the set of
nonnegative symmetric $d\times d$ real matrices.
The Fokker-Planck (or Kolmogorov forward) equation associated to \eqref{SDE} is
\begin{align}\label{PDE}
\partial_t f_t+ \hbox{div} (b(t,\cdot)f_t)=\frac 12 \sum_{i,j=1}^d \partial_{ij}(
[\sigma(t,\cdot) \sigma^*(t,\cdot)]_{i,j} f_t)
+ \ca L_t f_t,
\end{align}
where $\ca L_t f_t : \R^d\mapsto \R$ is defined by
$\int_{\R^d}(\ca L_t f_t)(x)\varphi(x) dx = \int_{\R^d}\int_E [\varphi(x+h(t,z,x))-\varphi(x)]f_t(x)dx$
for any reasonable $\varphi:\R^d\mapsto\R$. We use the notation $\nabla=\nabla_x$,
div$=$div$_x$ and $\partial_{ij}=\partial^2_{x_ix_j}$.

Let $\ca{P}(\R^d)$ be the set of probability measures on $\bb{R}^d$ and
\[\ca{P}_1(\R^d)=\{f\in\ca{P}(\R^d): m_1(f)<\infty\} \quad \text{with} \quad m_1(f):=\int_{\R^d} |x| f(dx).\]
We define  $L^\infty \big([0,T],\ca{P}_1(\R^d)\big)$ as the set of all measurable families
$(f_t)_{t\in[0,T]}$ of probability measures on $\R^d$ such that $\sup_{[0,T]}m_1(f_t)<\infty$.

\subsection{Main result}
We will suppose the following conditions.

\begin{assu}\label{Assum}
The functions
$\sigma : [0,T]\times \R^d\mapsto \ca{S}_d^+$, $b:[0,T]\times \R^d \mapsto \R^d$
and $h:[0,T]\times E \times \R^d  \mapsto \R^d$ are measurable and
there is a constant $C$ such that for
all $(t,x)\in [0,T]\times\R^d$,
\begin{align*}
|\sigma(t,x)| + |b(t,x)| + \int_{E} |h(t,z,x)|  \mu(dz)\le C (1+|x|).
\end{align*}
We set $a(t,x)=\sigma(t,x)\sigma^*(t,x)$, which satisfies $|a(t,x)|\leq C(1+|x|^2)$.
\end{assu}

\bdefi\label{ttt}
Suppose Assumption \ref{Assum}. A measurable family
$(f_t)_{t\in[0,T]}$ of probability measures on $\R^d$
is called a weak solution to \eqref{PDE} if for all  $\varphi\in C_c^2(\bb{R}^d)$, all $t\in [0,T]$,
\begin{align}\label{IPDE}
\int_{\bb{R}^d} \varphi(x)\, f_t(dx)=\int_{\R^d}\varphi(x)\, f_0(dx)
+\int_0^t\int_{\R^d} [ \ca{A}_s\varphi(x)+\ca{B}_s\varphi(x)]\, f_s(dx)\, ds,
\end{align}
with the diffusion operator
$\ca{A}_s\varphi(x):= b(s,x)\cdot \nabla \varphi(x)+\frac{1}{2}\sum_{i,j=1}^d a_{ij}(s,x)\partial_{ij}\varphi(x)$
and the jump operator $\ca{B}_s\varphi(x):=\int_{E}\big[\varphi(x+h(s,z,x))-\varphi(x)\big]\, \mu(dz).$
\edefi

We will check the following facts in the appendix, implying in particular that \eqref{IPDE} makes sense.

\begin{rem}\label{phiext} Suppose Assumption \ref{Assum}.

(i) For $\varphi\in C_c^2(\bb{R}^d)$, $\sup_{[0,T]\times\R^d}(|\ca{A}_s\varphi(x)|+|\ca{B}_s\varphi(x)|)<\infty$.

(ii) Any weak solution $(f_t)_{t\in[0,T]}$ to \eqref{PDE} starting from $f_0\in\ca P_1(\R^d)$
belongs to  $L^\infty([0,T],\ca{P}_1(\bb{R}^d))$.

(iii) If $f_0\in\ca P_1(\R^d)$,
the weak formulation \eqref{IPDE} automatically extends to all functions $\varphi \in C^2(\bb{R}^d)$
such that $(1+|x|)[|\varphi(x)|+|\nabla \varphi(x)|+|D^2\varphi(x)|]$ is bounded.
\end{rem}

Point (iii) is far from optimal, but sufficient for our purpose. Our main result reads as follows.

\bthm\label{thm-result}
Suppose Assumption \ref{Assum} and consider any weak solution $(f_t)_{t\in[0,T]}$
to \eqref{PDE} such that $f_0\in\ca P_1(\R^d)$.
There exist, on some probability space $(\Omega, \ca{F},(\ca{F}_t)_{t\in[0,T]},\bb{P})$,
a $d$-dimensional $(\ca{F}_t)_{t\in[0,T]}$-Brownian motion
$(B_t)_{t\in[0,T]}$, a \nc $(\ca{F}_t)_{t\in[0,T]}$-Poisson measure $N(dt,dz)$ on $[0,T]\times E$
with intensity measure $dt \,\mu(dz)$, these two objects being independent, as well as a
c\`{a}dl\`{a}g $(\ca{F}_t)_{t\in[0,T]}$-adapted process $(X_t)_{t\in[0,T]}$ solving
\eqref{SDE} and such that $\ca{L}(X_t)=f_t$ for all $t\in[0,T]$.
\ethm

For $(X_t)_{t\in [0,T]}$ a solution to \eqref{SDE} and for $f_t=\ca L(X_t)$, a simple application
of the It\^o formula (to compute $\int_{\R^d}\varphi(x)f_t(dx)=\E[\varphi(X_t)]$
with $\varphi\in C^2_c(\R^d)$) shows that the family $(f_t)_{t\in [0,T]}$ is a weak solution to \eqref{PDE}.
The following corollary is thus immediately deduced from Theorem \ref{thm-result}.

\bcor\label{ccc}
Suppose Assumption \ref{Assum} and fix $f_0\in\ca P_1(\R^d)$.

(i) The existence of a (weak) solution $(X_t)_{t\in[0,T]}$ to
\eqref{SDE} such that $\ca L(X_0)=f_0$ is equivalent to the existence 
of a weak solution $(f_t)_{t\in[0,T]}$
to \eqref{PDE} starting from $f_0$.

(ii) The uniqueness (in law) of the solution $(X_t)_{t\in[0,T]}$ to
\eqref{SDE}  with $\ca L(X_0)=f_0$ implies the uniqueness of the weak solution $(f_t)_{t\in[0,T]}$
to \eqref{PDE} starting from $f_0$.
\ecor

In almost all models arising from applied sciences, the jump operator is given under the form
$\ca B_s \varphi (x)= \int_F [\varphi(x+g(s,y,x)) - \varphi(x)]\kappa(s,y,x) \nu(dy)$,
meaning that when in the position $x$ at time $s$, the process jumps to $x+g(s,y,x)$
at rate $\kappa(s,y,x) \nu(dy)$. Here $F$ is a measurable space endowed with a $\sigma$-finite measure $\nu$
and we have two measurable functions $g : [0,T]\times F \times \R^d \mapsto \R^d$ and
$\kappa : [0,T]\times F \times \R^d \mapsto \R_+$. Introducing $E=F\times \R_+$, $\mu(dy,du)=\nu(dy)du$
and $h(s,(y,u),x)=g(s,y,x)\indiq_{\{u \leq \kappa(s,y,x)\}}$, one easily verifies  that
$\ca B_s \varphi (x)= \int_E [\varphi(x+h(s,(y,u),x)) - \varphi(x)] \mu(dy,du)$.
Our results thus apply if
$\int_F |g(s,y,x)|\kappa(s,y,x) \nu(dy) \leq C(1+|x|).$

\subsection{Motivation}
Stochastic differential equations with jumps are now playing an important role in modeling
and applied sciences. We refer to the book of Situ \cite{MR2160585} for all basic results
and a lot of possible applications.
The book of Jacod \cite{J} contains many important results about weak and strong existence and uniqueness,
relations between SDEs and martingale problems, etc. See also the survey paper of Bass \cite{B}.

Existence for PDEs is often more developed than for
SDEs, so Theorem \ref{thm-result} might be useful to derive some new weak existence results
for the SDE \eqref{SDE}.

Our main motivation is the uniqueness for some nonlinear PDEs, for which the use of nonlinear (in the sense
of McKean) SDEs has proved to be a powerful tool. For example, the first (partial)
uniqueness result concerning
the homogeneous Boltzmann for long range interactions was derived by Tanaka \cite{MR512334}.
He was studying the simplest case of Maxwell molecules. Unfortunately,
he was only able to prove the uniqueness in law of
the nonlinear SDE associated to the Boltzmann equation. Horowitz and Karandikar \cite{HK} were able
to deduce the uniqueness for the (same) Boltzmann equation proceeding as follows.
Let us recall that the original equation writes $\partial_t f_t =Q(f_t,f_t)$, for some quadratic nonlocal
operator $Q$.
For $f$ a solution, they consider the linear PDE $\partial_t g_t =Q(g_t,f_t)$, with unknown $g$
satisfying $g_0=f_0$.
They prove uniqueness in law for the (linear) SDE associated to this PDE (for any initial condition).
They deduce, extending some results of Ethier and Kurtz \cite[Chap.4, Propositions 9.18 and 9.19]{EK},
the uniqueness for the linear PDE (for any initial condition).
So the unique solution (with $g_0=f_0$) to $\partial_t g_t =Q(g_t,f_t)$ is $f$ itself.
Consequently,
the time marginals of the solution $X$ to the linear SDE (when $X_0\sim f_0$),
which solve $\partial_t g_t =Q(g_t,f_t)$ are necessarily $(f_t)_{t\in [0,T]}$. Thus $X$
actually solves the nonlinear SDE.
Since uniqueness in law holds for the nonlinear SDE by Tanaka \cite{MR512334},
they deduce that there is at most one solution to the Boltzmann equation $\partial_t f_t =Q(f_t,f_t)$,
for some given reasonable initial condition $f_0$.

Let us recall that the above mentioned results of Ethier and Kurtz (extended by Horowitz and
Karandikar \cite[Theorem B1]{HK} and by Bhatt and Karandikar \cite[e.g. Theorems 4.1 and 5.2 ]{BK})
state in spirit that if some SDE has a unique solution (in law) for any deterministic initial
condition, then the corresponding PDE has a unique weak solution for any reasonable initial condition.

Our result is much stronger, since it does not require at all  uniqueness  for \eqref{SDE}.
If, for example, studying the Boltzmann equation,
it directly implies that, to any solution $f$ to the nonlinear equation (seen here as a solution
to the linear equation $\partial_t g_t =Q(g_t,f_t)$), we can associate a solution
$X$ to the corresponding linear SDE with additionally $X_t\sim f_t$ for all $t$. In other words,
$X$ solves the nonlinear SDE. This might look anodyne, but this was crucial when studying
more singular nonlinear equations, such as the Landau or Boltzmann equations for moderately soft potentials,
see \cite{Fournier:2015aa} and \cite{Xu2016aa}. Indeed, in such cases, we really need to use some
physical symmetries to prove uniqueness : it is absolutely not clear that uniqueness holds for the linear
PDE $\partial_t g_t =Q(g_t,f_t)$, since one really uses that the two arguments of $Q$ are the same.

We hope the above discussion shows that Theorem \ref{thm-result}
is an interesting variation of the mentioned results
of Ethier and Kurtz \cite{EK}. As already said, the method we use
was initiated by Figalli \cite{MR2375067} for continuous SDEs ($h=0$) with bounded coefficients.
The boundedness assumption was relaxed in \cite[Appendix B]{Fournier:2015aa}.
A special jumping SDE (with $a=b=0$ and a special jump operator) was considered in \cite{Xu2016aa}
to study a singular homogeneous Boltzmann equation.
We decided to write down the general case in the present paper.
We did not want to assume some boundedness of the coefficients, although
it complicates the proofs without introducing new deep ideas, 
because it is very useful for practical purposes.

Finally, as explained in the next subsection, we are not able to prove a general result
when the jump part of the SDE has infinite variations, and this is a rather important limitation.

\subsection{Strategy of the proof and plan of the paper}

At many places, the situation is technically more
involved, but the global strategy is exactly the same as that introduced by Figalli \cite[Theorem 2.6]{MR2375067}.
Let  $(f_t)_{t\in[0,T]}$ be a given
weak solution to (\ref{PDE}).

I. In Section \ref{sec:reg}, we introduce $f_t^\e=f_t\star \phi_\e$, where $\phi_\e$ is
the centered Gaussian density with covariance matrix $\e I_d$.
We compute the PDE satisfied by $f_t^\e$: we find that
$\partial_t f_t^\e + \hbox{div} (b^\e(t,\cdot)f_t)=\frac 12 \sum_{i,j} \partial_{i,j}(a^\e_{i,j}(t,\cdot) f_t)+
\ca L_t^\e f_t^\e$,
for some coefficients $a^\e$, $b^\e$ and some jump operator $\ca L^\e_t$.
Let us mention that $a^\e(t,\cdot)$, $b^\e(t,\cdot)$ and $\ca L^\e_t$ of course depend on $f_t$.

II. Still in Section \ref{sec:reg}, we prove that $a^\e$, $b^\e$ and the coefficient of
the jump operator $\ca L^\e$ satisfy

(i)~the same linear growth conditions as $a$, $b$, $\ca L$, uniformly
in $\e\in (0,1)$,

(ii) some (non-uniform) local Lipschitz conditions.

III. In Section \ref{sec:eps}, we use II to build, for each $\e \in (0,1)$,
a solution $(X^\e_t)_{t\in [0,T]}$ to some SDE of which the  Fokker-Planck equation is the 
PDE satisfied by $(f_t^\e)_{t\in [0,T]}$. Since both the SDE and the PDE (with $\e \in (0,1)$ fixed)
are well-posed (because the coefficients are regular enough),
we conclude that $\ca L(X^\e_t)=f^\e_t$. Indeed, the time marginals of $(X^\e_t)_{t\in [0,T]}$
satisfy the same PDE as $(f_t^\e)_{t\in [0,T]}$.

IV. Still in Section \ref{sec:eps}, we prove that the family $\{(X^\e_t)_{t\in [0,T]},\e\in (0,1)\}$ is tight.
This is rather easy from the Aldous criterion \cite{MR0474446}, using only II-(ii).

V. In Section \ref{sec:conc}, we finally consider a limit point $(X_t)_{t\in [0,T]}$, as $\e\to 0$, of
$\{(X^\e_t)_{t\in [0,T]},\e\in (0,1)\}$. Since
$\ca L(X^\e_t)=f_t^\e$ by III, we deduce that $\ca L(X_t)=f_t$ for each $t\in [0,T]$.
It then remains to show that  $(X_t)_{t\in [0,T]}$ is a weak solution to \eqref{SDE} and
we classically make use of martingale problems.
Since the coefficients $a,b,h$ are possibly rough, we have to approximate them by some continuous (in $x$)
coefficients $\tilde a, \tilde b, \tilde h$. We use that we already know the time
marginals of $(X_t)_{t\in [0,T]}$: we can take $\tilde a(t,\cdot)$, $\tilde b(t,\cdot)$ and $\tilde h(t,\cdot,z)$
close to $a(t,\cdot)$, $b(t,\cdot)$ and $h(t,\cdot,z)$ in $L^1(f_t)$.

The proof of Remark \ref{phiext} is written in an appendix.

To conclude this paragraph, let us mention a few difficulties.
The regularized jump operator, in its weak form writes
$\int_{\R^d}\ca L_t^\e f_t^\e(y)\varphi(y)dy=
\int_{\R^d}\int_{\R^d}\int_E[\varphi(y+h(t,z,x))-\varphi(y)]\phi_\e(x-y)f_t^\e(dx)dy$.
We found no {\it regular} Poisson representation of the associated SDE. We use an indicator
function, see \eqref{convo-SDE}. This is why we are not able to treat the case of an
infinite variation jump term: we do not know how to prove that a SDE like \eqref{convo-SDE},
with a compensated Poisson measure and some weaker condition
on $h$ (something like $\int_E |h(s,z,x)|^2 \mu(dz) \leq C(1+|x|^2)$), is well-posed.

Although this should be classical since the coefficients are rather regular for $\e\in(0,1)$ fixed,
we found no reference about the uniqueness for the PDE satisfied by $(f_t^\e)_{t\in [0,T]}$
(see Lemma \ref{regu-eq}). We have not been able to write down a deterministic proof.
We thus use that the corresponding SDE is well-posed (for any deterministic initial condition) and
we apply a result of Horowitz and Karandikar \cite{HK}.

\subsection{Convention}

During the whole paper, we always suppose Assumption \ref{Assum} and that $f_0\in\ca P_1(\R^d)$.
We use the generic notation $C$ for a positive finite constant, of which the value 
may change from line to line. It is allowed to depend only on the dimension $d$, on the
parameters $a,b,h,E,\mu,T$ of our equations, and on the weak solution $(f_t)_{t\in [0,T]}$ to
\eqref{PDE} under study.
When a constant depends on another parameter, we indicate it in subscript.
For example, $C_\e$ is a constant allowed to depend only on
$a,b,h,E,\mu,T, (f_t)_{t\in [0,T]}$ and on $\e$.

\section{Regularization}\label{sec:reg}

We introduce the following regularization procedure, as Figalli in \cite{MR2375067},
see also \cite{Xu2016aa}.

\blem\label{regu-eq}
For $(f_t)_{t\in[0,T]} \in L^\infty([0,T],\ca P_1(\R^d))$ a weak solution to \eqref{PDE} and $\e\in(0,1)$,
we set
 \[f_t^\e(y):=\int_{\R^d}\phi_\e(x-y)f_t(dx)=(f_t\star\phi_\e)(y)
\quad with \quad \phi_\e(x)=(2\pi\e)^{-d/2}e^{-|x|^2/(2\e)}.\]
Then for any test function $\psi\in C_c^2(\R^d)$, any $t\in[0,T]$,
$$
\int_{\R^d}\psi(y)\,  f_t^\e(y)dy=\int_{\R^d}\psi(y)\,  f_0^\e(y)dy + \int_0^t \int_{\R^d}
[\ca{A}_{s,\e}\psi(y)+\ca{B}_{s,\e}\psi(y)] \, f_s^\e(y) dy ds,
$$
with
\begin{align*}
\ca{A}_{t,\e}\psi(y) =& b^\e(t,y)\cdot\nabla \psi(y)+\frac{1}{2}\sum_{i,j=1}^d a_{ij}^\e(t,y)\partial_{ij}\psi(y),\\
\ca{B}_{t,\e}\psi(y)=&\int_{E} \int_{\R^d}\big[\psi(y+h(t,z,x))-\psi(y)\big] F^\e_t(x,y) \,f_t(dx)\, \mu(dz),
\end{align*}
where
$$a^\e(t,y) =\frac{\int_{\R^d}\!\phi_\e(x-y)a(t,x) f_t(dx)}{f_t^\e(y)}, \; b^\e(t,y)
=\frac{\int_{\R^d}\!\phi_\e(x-y)b(t,x) f_t(dx)}{f_t^\e(y)},\;
F^\e_{t}(x,y):= \frac{\phi_\e(x-y)}{f_t^\e(y)}.
$$
\elem
\bpf
It is obvious that $f_t^\e(y)>0$ for each $(t,y)\in[0,T]\times \R^d$.
We first apply \eqref{IPDE} with the choice $\varphi(x)=\phi_\e(x-y)$ (with some fixed  $y\in \R^d$),
which is licit by Remark \ref{phiext}-(iii). We then integrate the obtained equality against
$\psi\in C_c^2(\R^d)$. This gives
$$
\int_{\R^d} \psi(y)f_t^\e(y)dy=\int_{\R^d} \psi(y)f_0^\e(y)dy + \int_0^t (I_s+J_s) ds,
$$
where
$$
I_t:=\int_{\R^d}\int_{\R^d} \psi(y)\ca{A}_t\phi_\e(x-y) f_t(dx) dy \quad \hbox{and}\quad
J_t:=\int_{\R^d}\int_{\R^d} \psi(y) \ca{B}_t\phi_\e(x-y) f_t(dx) dy.
$$
First,
$$
I_t= \int_{\R^d}\int_{\R^d}\psi(y) b(t,x)\cdot \nabla \phi_\e(x-y) f_t(dx) dy
+\frac{1}{2} \int_{\R^d}\int_{\R^d} \sum_{i,j=1}^d \psi(y)a_{ij}(t,x)\partial_{ij}\phi_\e(x-y) f_t(dx)dy.
$$
But we have $\int_{\R^d}\!\psi(y) \nabla \phi_\e(x-y) dy=\int_{\R^d} \phi_\e(x-y)  \nabla\psi(y)dy$
as well as $\int_{\R^d} \psi(y)\partial_{ij}\phi_\e(x-y) dy=\int_{\R^d} \phi_\e(x-y) \partial_{ij}\psi(y) dy$,
so that
\begin{align*}
I_t &=\int_{\R^d}\int_{\R^d} \phi_\e(x-y)b(t,x) \cdot \nabla\psi(y) f_t(dx) dy
+ \frac{1}{2} \int_{\R^d} \int_{\R^d} \sum_{i,j=1}^d a_{ij}(t,x) \phi_\e(x-y) \partial_{ij}\psi(y) f_t(dx) dy\\
&=\int_{\R^d} b^\e(t,y) \cdot \nabla\psi(y) f_t^\e(y) dy
+ \frac{1}{2} \int_{\R^d} \sum_{i,j=1}^d a_{ij}^\e(t,y)\partial_{ij}\psi(y) f_t^\e(y)dy\\
&=\int_{\R^d} \ca{A}_{t,\e}\psi(y) f_t^\e(y)dy
\end{align*}
as desired. For the jump term, we use a similar computation as in \cite[Proposition 3.1]{Xu2016aa}.
Since $\mu$ is $\sigma$-finite, there exists a non-decreasing sequence
$(E_n)_{n\ge1}\subset E$ such that $\bigcup_{n=1}^\infty E_n=E$ and $\mu(E_n)<\infty$ for each $n\ge1$.
We fix $n$ and write
\begin{align*}
J_t=&\int_{\R^d}\!\int_{\R^d}\!\int_{E_n} \psi(y)\phi_\e(x-y+h(t,z,x)) \mu(dz)f_t(dx)dy
- \int_{\R^d}\!\int_{\R^d}\!\int_{E_n} \psi(y) \phi_\e(x-y) \mu(dz)f_t(dx)dy\\
& +\int_{\R^d}\!\int_{\R^d}\!\int_{E\setminus E_n} \psi(y)\big[\phi_\e(x-y+h(t,z,x))-\phi_\e(x-y)\big]\mu(dz)f_t(dx)dy\,.
\end{align*}
Using the change of variables $y-h(t,z,x)\mapsto y$, we see that
$$
\int_{\R^d} \psi(y)\phi_\e(x-y+h(t,z,x)) dy
=\int_{\R^d}\psi(y+h(t,z,x))\phi_\e(x-y) dy,
$$
and consequently,
\begin{align*}
J_t=&\int_{\R^d}\int_{\R^d}\int_{E_n} \big[\psi(y+h(t,z,x))-\psi(y)\big]\phi_\e(x-y) \mu(dz) f_t(dx) dy\\
& +\int_{\R^d}\int_{\R^d}\int_{E\setminus E_n} \psi(y)\big[\phi_\e(x-y+h(t,z,x))-\phi_\e(x-y)\big]\mu(dz)f_t(dx) dy.
\end{align*}
Observe now that $|\psi(y+h(t,z,x))-\psi(y)|\phi_\e(x-y)
\leq C |h(t,z,x)|\phi_\e(x-y) \in L^1(\mu(dz) f_t(dx) dy)$
and $|\psi(y)[\phi_\e(x-y+h(t,z,x))-\phi_\e(x-y)]|\leq C_\e|\psi(y)||h(t,z,x)| \in L^1(\mu(dz) f_t(dx) dy)$:
this uses that $\psi \in C^2_c(\R^d)$, Assumption \ref{Assum} and
that $f_t\in \ca P_1(\R^d)$. We thus can let $n\to \infty$:
\begin{align*}
J_t=&\int_{\R^d}\int_{\R^d}\int_{E} \big[\psi(y+h(t,z,x))-\psi(y)\big]\phi_\e(x-y) \mu(dz) f_t(dx) dy\
=\int_{\R^d} \ca{B}_{t,\e}\psi(y) f_t^\e(y)dy,
\end{align*}
which completes the proof.
\epf

Let us now give some growth and regularity estimates on the regularized coefficients.

\blem\label{lem-condition}
Let  $(f_t)_{t\in[0,T]}\in L^\infty([0,T],\ca P_1(\R^d))$ be a weak solution to \eqref{PDE}
and recall that $a^\e, b^\e, F^\e$ were introduced in Lemma \ref{regu-eq}.

(i) There exists a constant $C>0$ such that for all $\e\in(0,1)$, all $y\in\bb{R}^d$, all $t\in[0,T]$,
\begin{align*}
|b^\e(t,y)|+|a^\e(t,y)|^{1/2}+
\int_{\R^d}  \int_E |h(t,z,x)| {F_t^\e(x,y)}   \mu(dz) f_t(dx)
\le C \,(1+|y|).
\end{align*}

(ii) For all $\e\in(0,1)$ and $R>0$, there is
$C_{R, \e}>0$  such that for all $y_1,y_2\in B(0, R)$, all $t\in[0,T]$,
\begin{align*}
|b^\e(t,y_1)-b^\e(t,y_2)| &+ |a^\e(t,y_1)-a^\e(t,y_2)|
+|[a^\e(t,y_1)]^{1/2}-[a^\e(t,y_2)]^{1/2}|\\
&+\int_{\R^d}  \int_{E}  |h(t,z,x)| |F_t^\e(x,y_1)- F_t^\e(x,y_2)|
\mu(dz) f_t(dx) \le C_{R, \e} \,|y_1-y_2|.
\end{align*}
\elem

\bpf
We start with (i). By Assumption \ref{Assum},
\begin{align*}
&|b^\e(t,y)|+|a^\e(t,y)|^{1/2}+
\int_{\R^d}  \int_E |h(t,z,x)| F_t^\e(x,y)   \mu(dz) f_t(dx) \\
\leq &  C  \frac{\int_{\R^d}\phi_\e(x-y)(1+|x|)\, f_t(dx)}{f_t^\e(y)}
+ C \Big[\frac{\int_{\R^d}\phi_\e(x-y)(1+|x|)^2\, f_t(dx)}{f_t^\e(y)}\Big]^{1/2}\\
=: &CI_\e(t,y)+CJ_\e(t,y).
\end{align*}
Since for $y$ fixed, $[f_t^\e(y)]^{-1}\phi_\e(x-y)f_t(dx)$ is a probability measure, we infer from
Cauchy-Schwarz that $I_\e(t,y)\leq J_\e(t,y)$. We thus only have to prove that
$[J_\e(t,y)]^2 \leq C(1+|y|^2)$.
Let $L:=2\sup_{[0,T]}m_1(f_t)+2$.
We use that
$$
1+|x| \leq 1+|y|+|x-y|\leq 1+2|y|+L+|x-y| \indiq_{\{|x-y|>|y|+L\}}
$$
to write
\begin{align*}
[J_\e(t,y)]^2 \le &
2 \frac{\int_{\R^d}(1+2|y|+L)^2\phi_\e(x-y) f_t(dx)}{f_t^\e(y)}
+2 \frac{\int_{|x-y|\ge |y|+L} |x-y|^2\phi_\e(x-y) f_t(dx)}{f_t^\e(y)}\\
\leq &2(1+2|y|+L)^2 +2 \frac{(|y|+L)^2\phi_\e(|y|+L)}{f_t^\e(y)}.
 \end{align*}
For the second term, we used that $|y|+L\geq 2 \geq \sqrt{2\e}$ and that
$z \mapsto |z|^2 \,\phi_\e(z)$ is radially symmetric and
decreasing on $\{|z|\geq \sqrt{2\e}\}$. To conclude the proof of (i), it suffices to note
that
$$
f_t^\e(y)\ge\int_{|x-y|\le |y|+L} \phi_\e(x-y) \,f_t(dx) \ge \phi_\e(|y|+L) \,f_t(B(y, |y|+L))
\ge \phi_\e(|y|+L)/2
$$
because $z \mapsto \phi_\e(z)$ is radially symmetric decreasing and because
$f_t(B(y, |y|+L))\geq f_t(B(0,L)) \geq 1/2$, since $f_t(B(0,L)^c)\leq m_1(f_t)/L \leq 1/2$.

For point (ii), it suffices to prove that
$\nabla_y b^\e(t,y)$, $\nabla_y a^\e(t,y)$, $D^2_y a^\e(t,y)$ are locally bounded on $[0,T]\times \R^d$, as well as
$G^\e(t,y):=\int_{\R^d}\int_E |h(t,z,x)||\nabla_y F_{t}^\e(x,y)|\mu(dz) f_t(dx)$.
No uniformity in $\e$ is required here. By Stroock and Varadhan \cite[Theorem 5.2.3]{SW},
the local boundedness of $D^2_y a^\e(t,y)$ implies that of $\nabla_y ([a^\e(t,y)]^{1/2})$.

First, one easily checks that $y\mapsto (f_t^\e(y))^{-1}$ is of class $C^\infty$ for each $t\in [0,T]$
and that it is locally bounded, as well as its derivatives of order $1$ and $2$, on $[0,T]\times\R^d$.
This uses in particular the lower bound $f_t^\e(y)\ge \phi_\e(|y|+L)/2$ proved a few lines above.

Recall that by definition, we have
$a^\e(t,y)=(f_t^\e(y))^{-1}\int_{\R^d} \phi_\e(x-y)a(t,x)f_t(dx)$ and
$b^\e(t,y)=(f_t^\e(y))^{-1}\int_{\R^d} \phi_\e(x-y)b(t,x)f_t(dx)$.
Recall finally that $|a(t,x)|+|b(t,x)|\leq C(1+|x|^2)$.
So concerning $a^\e$ and $b^\e$, our goal is only to check that
$$
K_\e(t,y):=\int_{\R^d} [|\nabla_y \phi_\e(x-y)|+ |D^2_y \phi_\e(x-y)|](1+|x|^2) f_t(dx)
$$
is locally bounded on $[0,T]\times \R^d$.
But using that $ (1+|z|^2)  [|\nabla \phi_\e(z)|+ |D^2 \phi_\e(z)|]$ is bounded on $\R^d$, we
deduce that $[|\nabla_y \phi_\e(x-y)|+ |D^2_y \phi_\e(x-y)|](1+|x|^2)\leq C_\e(1+|y|^2)$,
whence $K_\e(t,y) \leq C_\e(1+|y|^2)$.

Next, one has $|\nabla_y F_t^\e(x,y)|\leq C_\e (f_t^\e(y))^{-2}[\phi_\e(x-y) |\nabla f^\e_t(y)|
+ f^\e_t(y)|\nabla \phi_\e(x-y)|]$. Using again that $f_t^\e$ is smooth and positive, the goal
concerning $G^\e$ is to verify that
$$
L_\e(t,y):= \int_{\R^d}\int_E |h(t,z,x)|[\phi_\e(x-y) + |\nabla \phi_\e(x-y)|]\mu(dz) f_t(dx)
$$
is locally bounded. By Assumption \ref{Assum},
$$
L_\e(t,y) \leq \int_{\R^d}[\phi_\e(x-y) + |\nabla \phi_\e(x-y)|](1+|x|) f_t(dx) \leq C_\e(1+|y|)
$$
as previously, because $(1+|z|)[\phi_\e(z) + |\nabla \phi_\e(z)|]$ is  bounded.
\epf

\section{Study of the regularized equations}\label{sec:eps}

In this section, we build a realization of the regularized weak solution $(f_t^\e)_{t\in[0,T]}$.

\bprop\label{well-posedness}
Let $(f_t)_{t\in[0,T]}\in L^\infty([0,T],\ca P_1(\R^d))$ be a weak  solution to \eqref{PDE} and fix $\e\in (0,1)$.
Consider $(f_t^\e)_{t\in[0,T]}$
and  $a^\e, b^\e, F^\e$ defined in Lemma \ref{regu-eq} and put
$\sigma^\e(t,y):=(a^\e(t,y))^{1/2}$.
Consider a random variable $X_0^\e$, a $d$-dimensional Brownian motion
$(B_s)_{s\in[0,T]}$ and a Poisson measure $N(ds, dz, dx, du)$
on $[0,T]\times E\times \R^d\times[0,\infty)$ with intensity measure $ds\, \mu(dz)\, f_s(dx)\, du$,
these three objects being independent. We work with the filtration generated by
$X_0^\e,B,N$.

(i) There is a pathwise unique c\`adl\`ag adapted solution $(X_t^\e)_{t\in[0,T]}$ to
\begin{align}\label{convo-SDE}
X_t^\e=&X_0^\e+\int_0^t b^\e(s,X_s^\e) ds + \int_0^t \sigma^\e(s,X_s^\e) dB_s \notag\\
&+ \int_0^t\int_E\int_{\R^d}\int_0^\infty h(s,z,x) \bbd{1}_{\{u\le F_s^\e(x, X_{s-}^\e)\}}N(ds, dz, dx, du).
 \end{align}

(ii) There is a constant $C$ (not depending on $\e$) such that $\bb{E}[\sup_{[0,T]}|X_t^\e|]
\leq C(1+\E[|X_0^\e|])$.

(iii) If $\ca L(X_0^\e)=f_0^\e$, then $\ca{L}(X_t^\e)=f_t^\e$ for all $t\in[0,T]$.
\eprop
\bpf
(i) The existence of a pathwise unique solution to \eqref{convo-SDE} is more or less standard,
because of the linear growth and local Lipschitz properties of the coefficients
proved in Lemma \ref{lem-condition}. We only prove pathwise uniqueness, the existence
being shown similarly, using a localization procedure (to make the coefficients
globally Lipschitz continuous) and a Picard iteration.
Consider two solutions $(X_t^\e)_{t\in[0,T]}$ and $(\tilde{X}_t^\e)_{t\in[0,T]}$ to
\eqref{convo-SDE} with $X_0^\e=\tilde X_0^\e$ and introduce the stopping time
$\tau_R:=\inf\{t\in[0,T] : |X_t^\e|\vee |\tilde{X}_t^\e|\ge R\}$, for $R>0$, with the convention
that $\inf \emptyset = T$. Using the Burkholder-Davis-Gundy inequality for the Brownian
part, we find
\begin{align*}
\E\Big[\sup_{[0,t\land \tau_R]} |X^\e_s-\tilde X^\e_s|\Big]\leq & \E\Big[\int_0^{t\land\tau_R}\!\!\!\! |b^\e(s,X_s^\e)
-b^\e(s,\tilde X_s^\e)| ds + C\Big(\int_0^{t\land\tau_R}\!\!\!\! |\sigma^\e(s,X_s^\e)-\sigma^\e(s,\tilde X_s^\e)|^2
ds \Big)^{1/2} \\
&\hskip0.8cm + \int_0^{t\land\tau_R}\! \int_E \int_{\R^d}|h(s,z,x)| |F_s^\e(x, X_{s}^\e)- F_s^\e(x, \tilde X_{s}^\e)|
f_s(dx) \mu(dz) ds \Big].
\end{align*}
By Lemma \ref{lem-condition}-(ii), we deduce that
\begin{align*}
\E\Big[\sup_{[0,t\land \tau_R]} |X^\e_s-\tilde X^\e_s|\Big]\leq & C_{R,\e}\E\Big[\int_0^{t\land\tau_R} |X_s^\e
-\tilde X_s^\e| ds + \Big(\int_0^{t\land\tau_R}\! |X_s^\e-\tilde X_s^\e|^2 ds \Big)^{1/2}\Big]\\
\leq & C_{R,\e} (t+\sqrt t) \E\Big[\sup_{[0,t\land \tau_R]} |X^\e_s-\tilde X^\e_s|\Big].
\end{align*}
We deduce that $\E[\sup_{[0,t_R\land \tau_R]} |X^\e_s-\tilde X^\e_s|]=0$, where $t_R>0$ is
such that $C_{R,\e} (t_R+\sqrt {t_R})=1/2$.
But then, the same computation allows us to prove that
$\E[\sup_{[t_R\land \tau_R,(2t_R)\land \tau_R]} |X^\e_s-\tilde X^\e_s|]=0$, etc, so that we end with
$\E[\sup_{[0,T\land \tau_R]} |X^\e_s-\tilde X^\e_s|]=0$ for each $R>0$.
Since $\lim_{R\to\infty}\tau_R=T$ a.s. (because $(X_t^\e)_{t\in[0,T]}$ and $(\tilde{X}_t^\e)_{t\in[0,T]}$
are assumed to be a.s. c\`adl\`ag and thus locally bounded on $[0,T]$),
we conclude that $\E[\sup_{[0,T]} |X^\e_s-\tilde X^\e_s|]=0$, which was our goal.

(ii)
Using the Burkholder-Davis-Gundy inequality for the Brownian part, we find, for $t\in[0,T]$,
\begin{align*}
u_t^\e:=\bb{E}\Big[\sup_{[0,t]}|X_s^{\e}|\Big] \le& \bb{E}[|X_0^\e|]
+ \bb{E}\Big[\int_0^t |b^\e(s, X_s^{\e})| ds\Big] + C \bb{E}\Big[\Big(\int_0^t |\sigma^\e(s,X_s^{\e})|^2 \,
ds\Big)^{1/2}\Big]\\
&+ \bb{E}\Big[\int_0^t\int_{E}\int_{\R^d}|h(s,z,x)|F_{s}^\e(x, X_{s-}^{\e}) f_s(dx)\mu(dz)ds\Big].
\end{align*}
Inserting the estimates proved in Lemma \ref{lem-condition}-(i), we find, for some
constant $C$ not depending on $\e\in(0,1)$ nor on $\bb{E}[|X_0^\e|]$,
$$
u_t^\e \le \bb{E}[|X_0^\e|]+C\bb{E}\Big[\int_0^t \big(1+|X_s^{\e}| \big)ds + \Big(\int_0^t
(1+|X_s^{\e}|^2) \, ds\Big)^{1/2}\Big] \leq u_0^\e  +   C(t+\sqrt t) (1+u_t^\e).
$$
With $t_0>0$ such that $C(t_0+\sqrt {t_0})=1/2$, we conclude that $u_{t_0}^\e \leq 2u_0^\e +1$.
One checks similarly that $u_{2t_0}^\e \leq 2u_{t_0}^\e +1\leq 4u_0^\e +3$.
Repeating the argument, we end with
$u_T^\e \leq 2^{\lfloor T/t_0\rfloor+1} u_0^\e +2^{\lfloor T/t_0\rfloor+1}-1$.

(iii) We now assume that $\ca L(X_0^\e)=f_0^\e$ and we set $g_t^\e:=\ca{L}(X_t^\e)$.
A direct application of the It\^o formula shows that for all $t\in[0,T]$,
recalling the notation of Lemma \ref{regu-eq},
\begin{align*}
\int_{\R^d}\psi(y)\,  g_t^\e(dy)=&\int_{\R^d}\psi(y)\,  f_0^\e(dy) +
\int_0^t \int_{\R^d} [\ca{A}_{s,\e}\psi(y)+\ca{B}_{s,\e}\psi(y)] \, g_s^\e(dy) ds.
\end{align*}
Recalling  Lemma \ref{regu-eq} again, $(f_t^\e)_{t\in[0,T]}$ solves the same equation.
The following uniqueness result will thus complete the proof of (iii): for any
$\nu_0\in\ca{P}(\R^d)$, there exists at most one measurable family
$(\nu_t)_{t\in [0,T]}$ of probability measures such that for all $\psi\in C_c^2(\R^d)$ and all $t\in [0,T]$,
\be{martingale-problem}
\int_{\R^d}\psi(y)\, \nu_t(dy)=\int_{\R^d}\psi(y)\, \nu_0(dy)+\int_0^t ds\int_{\R^d} \nu_s(dy)\,
[\ca{A}_{s,\e}\psi(y)+\ca{B}_{s,\e}\psi(y)]\,.
\ee
This must be classical (because the coefficients are rather regular),
but we found no reference and thus make use of martingale problems. A c\`{a}dl\`{a}g adapted
$\R^d$-valued process $(Y_t)_{t\in[0,T]}$ on some filtered probability space
$(\Omega, \ca{F},(\ca{F}_t)_{t\in[0,T]},\bb{P})$
is said to solve $MP_\e(\nu_0)$ if
$\ca L(Y_0)=\nu_0$ and if
$$
\psi(Y_t)-\int_0^t \big[\ca{A}_{s,\e}\psi(Y_s)+\ca{B}_{s,\e}\psi(Y_s)\big]\, ds
$$
is a martingale for all $\psi\in C_c^2(\R^d)$.
Due to Horowitz and Karandikar \cite[Theorem B1]{HK}, the following points
imply uniqueness for \eqref{martingale-problem}.
Here $C_0(\R^d)$ is the set of continuous functions from $\R^d$ to $\R$ vanishing at infinity.

(a) $C^2_c(\R^d)$ is dense is $C_0(\R^d)$ for the uniform convergence topology,

(b) $(t,y)\mapsto \ca A_{t,\e}\psi(y)+\ca B_{t,\e}\psi(y)$ is measurable for all $\psi\in C^2_c(\R^d)$,

(c) for each $t\in[0,T]$, $\ca A_{t,\e}+\ca B_{t,\e}$ satisfies the maximum principle,

(d) there exists a countable family $(\psi_k)_{k\ge1}\subset C_c^2(\R^d)$ such that for all $t\in[0,T]$,
\[\overline{\{(\psi_k,  \ca{A}_{t,\e}\psi_k+\ca{B}_{t,\e}\psi_k), k\ge 1\}}  ~\supset~
\{(\psi,  \ca{A}_{t,\e}\psi+\ca{B}_{t,\e}\psi), \psi\in C_c^2(\R^d)\}\] where  the closure in the left-hand side
is under the bounded pointwise convergence,

(e) for each $y_0\in\R^d$, there exists a unique (in law) solution to $MP_\e(\delta_{y_0})$.

Points (a) and (b) are obvious. The SDE associated to $MP_\e$ is precisely \eqref{convo-SDE}:
$(Y_t)_{t\in[0,T]}$ solves $MP_\e(\nu_0)$ if and only if it is a weak solution to \eqref{convo-SDE}
and $\ca L(Y_0)=\nu_0$, see Jacod \cite[Theorem 13.55]{J}, see also \cite[Theorem A1]{HK}.
Thus (e) follows from (i).
For (c), assume that $\psi \in C^2_c(\R^d)$ attains its maximum at $y_0$.
Then $\ca B_{t,\e} \psi(y_0)\leq 0$ (this is immediate) and $\ca A_{t,\e} \psi(y_0)\leq 0$
(because $\nabla \psi(y_0)=0$ and, since $a(t,y_0)$ is symmetry and nonnegative,
$\sum_{i,j} a_{ij}(t,y_0)\partial_{ij}\psi(y_0)\leq 0$). It only remains to prove (d).
Consider any countable subset
$(\psi_k)_{k\ge1}\subset C_c^2(\R^d)$ dense in $C_c^2(\R^d)$: for $\psi\in C_c^2(\R^d)$ with
Supp $\psi\subset B(0,M)$,
there exists $(\psi_{k_n})_{n\ge1}$ with Supp $\psi_{k_n}\subset B(0,2M)$ such that
\[\lim_{n\to\infty}(\|\psi-\psi_{k_n}\|_\infty+\|\nabla (\psi-\psi_{k_n})\|_\infty+\|D^2(\psi-\psi_{k_n})\|_\infty)=0\,.\]
We will prove more than needed, namely that
(i) $\lim_{n\to\infty} \sup_{[0,T]}\|\ca{A}_{t,\e}\psi_{k_n}-\ca{A}_{t,\e}\psi\|_\infty=0$,
and (ii) $\lim_{n\to\infty}\sup_{[0,T]} \|\ca{B}_{t,\e}\psi_{k_n}-\ca{B}_{t,\e}\psi\|_\infty=0$.

By Lemma \ref{lem-condition},
\begin{align*}
|\ca{A}_{t,\e}(\psi_{k_n}-\psi)(y)|
&\le \|\nabla (\psi_{k_n}-\psi)\|_\infty \,|b^\e(t,y)| \,\bbd{1}_{\{|y|\le 2M\}} + \frac{1}{2}\,
\|D^2(\psi_{k_n}-\psi)\|_\infty \,\|a^\e(t,y)\| \,\bbd{1}_{\{|y|\le 2M\}} \\
&\leq C \|\nabla(\psi_{k_n}-\psi)\|_\infty +C\|D^2(\psi_{k_n}-\psi)\|_\infty,
\end{align*}
which tends to $0$, implying (i). We next write, using that Supp $(\psi_{k_n}-\psi)\subset B(0,2M)$,
\begin{align*}
|(\psi_{k_n}-\psi)(y+h(t,z,x))-(\psi_{k_n}-\psi)(y)|\leq& \bbd{1}_{\{|y|\le 4M\}}
\|\nabla(\psi_{k_n}-\psi)\|_\infty|h(t,z,x)|\\
&+ 2\bbd{1}_{\{|y|\ge 4M\}} \|\psi_{k_n}-\psi\|_\infty\bbd{1}_{\{|y+h(t,x,z)|\le 2M\}}.
\end{align*}
Observing that
\[\bbd{1}_{\{|y|\ge 4M, |y+h(t,z,x)|\le 2M\}}\le \bbd{1}_{\{|y|\ge 4M, |h(t,z,x)|\ge |y|/2\}} \le \bbd{1}_{\{|y|\ge 4M\}}
\,\frac{2|h(t,z,x)|}{|y|},\]
we deduce that
\begin{align*}
|\ca{B}_{t,\e}(\psi_{k_n}-\psi)(y)|\le & \,\bbd{1}_{\{|y|\le 4M\}}\, \|\nabla(\psi_{k_n}-\psi)\|_\infty
\int_E\int_{\R^d}|h(t,z,x)|F_t^\e(x,y)\, f_t(dx)\, \mu(dz)\\
&+\bbd{1}_{\{|y|\ge 4M\}} \|\psi_{k_n}-\psi\|_\infty\int_E\int_{\R^d} \frac{2|h(t,z,x)|}{|y|}
F_t^\e(x,y) \,f_t(dx)\,\mu(dz).
\end{align*}
Recalling that $\int_E\int_{\R^d}|h(t,z,x)| F_t^\e(x,y) \,f_t(dx)\,\mu(dz)\leq C(1+|y|)$ by Lemma
\ref{lem-condition}, we find
\begin{align*}
|\ca{B}_{t,\e}(\psi_{k_n}-\psi)(y)|\le & \,\bbd{1}_{\{|y|\le 4M\}}C(1+|y|) \|\nabla(\psi_{k_n}-\psi)\|_\infty
+ \bbd{1}_{\{|y|\ge 4M\}} C \frac{1+|y|}{|y|} \|\psi_{k_n}-\psi\|_\infty\\
\leq& C  \|\nabla(\psi_{k_n}-\psi)\|_\infty + C  \|\psi_{k_n}-\psi\|_\infty
\end{align*}
and the conclusion follows.
\epf

\blem\label{lem-tightness}
For $(f_t)_{t\in[0,T]}\in L^\infty([0,T],\ca P_1(\R^d))$ a weak  solution to \eqref{PDE}
and $\e\in(0,1)$, consider the process $(X_t^\e)_{t\in [0,T]}$, with $X_0^\e\sim f_0^\e$, introduced in Lemma
\ref{well-posedness}. The family $\{(X_t^\e)_{t\in [0,T]},\e>0\}$ is tight in
$\bb{D}([0,T],\R^d)$ and any limit point $(X_t)_{t\in [0,T]}$ satisfies
$\bb{P}(\Delta X_t \neq 0)=0$ for all $t\in[0,T]$.
\elem

\bpf
We use the Aldous criterion \cite{MR0474446}, see also Jacod and Shiryaev \cite[p. 356]{MR1943877},
which implies tightness and that any limit point $(X_t)_{t\in [0,T]}$
is quasi-left-continuous and thus has no deterministic jump time.
It suffices to check that

(i) $\sup_{\e\in(0,1)}\bb{E}[\sup_{[0,T]}|X_t^{\e}|]<\infty$,

(ii) $\lim_{\beta\to0}~ \sup_{\e\in(0,1)}\sup_{(S,S^\prime)\in\ca{S}_T(\beta)}~ \bb{E}[|X_{S^\prime}^{\e}-X_{S}^{\e}|]=0$,
 where $\ca{S}_T(\beta)$ is the set of all pairs of stopping times $(S,S^\prime)$ satisfying
$0\le S\le S^\prime\le S+\beta\le T$ a.s.

Point (i) has already been checked in Lemma \ref{well-posedness}-(ii), since
$\E[|X_0^\e|]=m_1(f_0^\e) \leq m_1(f_0)+\sqrt{d\e}$.
Next, for $S,S^\prime\in\ca{S}_T(\beta)$ and $\e\in (0,1)$, we have
\begin{align*}
\bb{E}[|X_{S^\prime}^{\e}-X_{S}^{\e}|]\le
& \bb{E}\Big[\int_S^{S+\beta} |b^\e(s,X_s^{\e})| \,ds\Big] + \bb{E}\Big[\Big|\int_S^{S^\prime}
\sigma^\e(s,X_{s}^{\e}) \,dB_s\Big| \Big] \\
& + \bb{E}\Big[\int_S^{S+\beta}\int_{E}\int_{\R^d}|h(s,z,x)|F_s^\e(x, X_{s}^{\e}) \,f_s(dx)\, \mu(dz)\, ds\Big]\\
\le& C\, \bb{E}\Big[\int_S^{S+\beta} \big(1+|X_s^\e|\big)\, ds\Big] +
C\bb{E}\Big[\Big(\int_S^{S^\prime} |\sigma^\e(s,X_{s}^{\e})|^2 ds \Big)^{1/2} \Big],
\end{align*}
where the last inequality follows from Lemma \ref{lem-condition}-(i) and the
Burkholder-Davis-Gundy inequality. But $|\sigma^\e(s,x)|^2 \leq C|a^\e(s,x)|\leq C(1+|x|^2)$ by
Lemma \ref{lem-condition}-(i) again, whence
$$
\bb{E}[|X_{S^\prime}^{\e}-X_{S}^{\e}|]
\le C\, \bb{E}\Big[\int_S^{S+\beta} (1+|X_s^\e|)\, ds +
\Big(\int_S^{S+\beta} (1+|X_s^\e|^2) ds \Big)^{1/2} \Big].
$$
Hence
$\bb{E}[|X_{S^\prime}^{\e}-X_{S}^{\e}|]
\leq C(\beta+\sqrt\beta)\bb{E}[\sup_{[0,T]}(1+|X_s^\e|)] \leq C (\beta+\sqrt\beta)$,
which ends the proof.
\epf

\section{Conclusion}\label{sec:conc}

As Figalli \cite{MR2375067}, we will need some continuous (in $x$) approximations of $a$, $b$ and $h$.

\blem\label{lem-continuous}
Let $(f_t)_{t\in[0,T]}\in L^\infty([0,T],\ca P_1(\R^d))$ be a weak  solution to \eqref{PDE}.
For all $\rho>0$, we can find  $\tilde a :[0,T]\times\R^d\mapsto \ca S_d^+$
and $\tilde b :[0,T]\times\R^d\mapsto \R^d$, both continuous and compactly supported,
a set $A\in\ca E$ such that $\mu(A)<\infty$, and a measurable
function $\tilde h : [0,T]\times E\times\R^d\mapsto \R^d$, continuous on $[0,T]\times \R^d$ for each $z\in E$,
such that $\tilde h(t,z,x)=0$ for all $(t,z,x)\in [0,T]\times A^c \times \R^d$ and
$$
\int_0^T \int_{\R^d} \Big[\frac{|a(t,x)-\tilde a(t,x)|}{1+|x|} + |b(t,x)-\tilde b(t,x)| +
\int_E |h(t,z,x)-\tilde h(t,z,x)|\mu(dz)
\Big] f_t(dx)dt < \rho.
$$
\elem

\bpf For $a$ and $b$, this follows from the fact, see Rudin \cite[Theorem 3.14]{R},
that continuous functions with compact support
are dense in $L^1([0,T]\times\R^d, dt f_t(dx))$, and that both $a(t,x)/(1+|x|)$ and $b(t,x)$
belong to this space by Assumption \ref{Assum}.

Since $h \in L^1([0,T]\times E\times\R^d,dt \mu(dz) f_t(dx))$ by Assumption \ref{Assum}
and since $\mu$ is $\sigma$-finite, we can find $A\in \ca E$ such that $\mu(A)<\infty$ and
$\int_0^T\int_{A^c}\int_{\R^d} |h(t,z,x)| f_t(dx)\mu(dz) dt < \rho/3$.

Next, can find a simple function $g=\sum_{n=1}^N \alpha_n \indiq_{S_n}$, with $\alpha_n \in \R_*$,
$S_n\in \ca B([0,T]\times \R^d)\otimes\ca E$, such that
$\int_0^T\int_A\int_{\R^d} |g(t,z,x)- h(t,z,x)| f_t(dx)\mu(dz) dt < \rho/3$.

But for $S\in \ca B([0,T]\times \R^d)\otimes\ca E$ and
$\e>0$, there is $\varphi_{S,\e}:[0,T]\times \R^d \times E \mapsto \R$, measurable,
continuous on $[0,T]\times \R^d$ for each $z\in E$ and
such that $\int_0^T\int_A\int_{\R^d} |\indiq_{\{(t,z,x)\in S\}}-
\varphi_{S,\e}(t,z,x)| f_t(dx)\mu(dz) dt < \e$.
Indeed, when $S=C\times D$ with $C\in\ca B([0,T]\times \R^d)$ and $D\in \ca E$,
it suffices to consider $\psi$ continuous on $[0,T]\times \R^d$ such that
$\int_0^T\int_{\R^d} |\indiq_{\{(t,x)\in C\}}-\psi(t,x)| f_t(dx)dt < \e/\mu(A)$ and to set
$\varphi_{S,\e}(t,z,x)= \psi(t,x)\indiq_{\{z \in D\}}$. The general case follows from the monotone 
class theorem.

Finally, $\tilde h (t,z,x)= \sum_{n=1}^N \alpha_n \varphi_{S_n,\rho/(3 |\alpha_n| 2^n)}(t,z,x)\indiq_{\{z\in A\}}$
is measurable and continuous in $(t,x)$ for each $z\in E$. Writing
\begin{align*}
|h(t,z,x)-\tilde h(t,z,x)|\leq& |h(t,z,x)|\indiq_{\{z\in A^c\}} +  |g(t,z,x)- h(t,z,x)|\indiq_{\{z\in A\}}\\
&+ \sum_{n=1}^N |\alpha_n||\varphi_{S_n,\rho/(3 |\alpha_n| 2^n)}(t,z,x) -  \indiq_{\{(t,z,x)\in S_n\}}|\indiq_{\{z\in A\}},
\end{align*}
we conclude that
$\int_0^T \int_E\int_{\R^d} |h(t,z,x)-\tilde h(t,z,x)| f_t(dx)\mu(dz) dt <\rho$
as desired.
\epf

We now can give the

\bpf[Proof of Theorem \ref{thm-result}]
Let  $(f_t)_{t\in[0,T]}\in L^\infty([0,T],\ca P_1(\R^d))$ be a weak solution to \eqref{PDE}.
For each $\e\in(0,1)$, consider $(f^\e_t)_{t\in[0,T]}$ introduced in Lemma \ref{regu-eq} and
the process $(X^\e_t)_{t\in[0,T]}$ introduced in Lemma \ref{well-posedness}-(iii).
By Lemma \ref{lem-tightness}, we can find a sequence $(X_t^{\e_n})_{t\in [0,T]}$
converging in law to some
process $(X_t)_{t\in[0,T]}$. Since we know from Lemma \ref{well-posedness} that $\ca{L}(X_t^{\e_n})=f_t^{\e_n}$
for each $t\in[0,T]$, each $n\geq 1$ and since
$f_t^{\e_n}$ goes weakly to $f_t$ as $n\to \infty$ by construction, we deduce that
for all $t\in [0,T]$, $\ca L(X_t)=f_t$. It thus only remains to verify that $X:=(X_t)_{t\in[0,T]}$
is a (weak) solution to \eqref{SDE}. According to the theory of martingale problems,
see Jacod \cite[Theorem 13.55]{J},
it classically suffices to prove that for any $\psi\in C_c^2(\R^d)$, the process
\[\psi(X_t)-\psi(X_0)-\int_0^t \big[ \ca{A}_s\psi(X_s)+\ca{B}_s\psi(X_s)\big]\, ds\]
is a martingale in the filtration $\ca F_t=\sigma(X_s,s\leq t)$. Our goal is thus to check that
for any $0\le s_1\le \dots \le s_k \le s\le t\le T$, any $\psi_1,\dots, \psi_k\in C_b(\R^d)$
and any $\psi\in C_c^2(\R^d)$, we have
$\bb{E}[\ca{K}(X)]=0$,
where $\ca{K}: \bb{D}([0,T], \R^d)\mapsto \R$ is defined by
\[\ca{K}(\lambda):=\Big(\prod_{i=1}^{k}\psi_i(\lambda_{s_i})\Big) \Big(\psi(\lambda_t)-\psi(\lambda_s)
-\int_s^t \big[\ca{A}_r\psi(\lambda_r)+\ca{B}_r\psi(\lambda_r)\big] \,dr\Big)\,.\]

We fix $\rho>0$ and consider $\tilde{a}$, $\tilde{b}$ and $\tilde h$ introduced in Lemma \ref{lem-continuous}.
We introduce $\tilde{\ca{A}}_s$ and $\tilde{\ca{B}}_s$ exactly as in Definition \ref{ttt} with
$\tilde{a}$, $\tilde{b}$ and $\tilde h$ instead of $a$, $b$ and $h$.
We define $\tilde{a}^\e$, $\tilde{b}^\e$, $\tilde{\ca{A}}_{s,\e}$ and $\tilde{\ca{B}}_{s,\e}$
exactly as in Lemma \ref{regu-eq}, with everywhere $\tilde{a}$, $\tilde{b}$ and $\tilde h$
instead of $a$, $b$ and $h$. Finally, we define $\tilde{\ca{K}}$
(resp. $\tilde{\ca{K}}_\e$, resp. $\ca K_\e$)
exactly as $\ca K$ with $\ca{A}_r$ and $\ca{B}_r$ replaced by $\tilde{\ca{A}}_r$ and $\tilde{\ca{B}}_r$
(resp. by  $\tilde{\ca{A}}_{r,\e}$ and $\tilde{\ca{B}}_{r,\e}$, resp. by ${\ca{A}}_{r,\e}$ and ${\ca{B}}_{r,\e}$).

First, $\bb{E}[\ca{K}_{\e_n}(X^{\e_n})]=0$. Indeed, since $X^{\e}=(X_t^{\e})_{t\in [0,T]}$
solves \eqref{convo-SDE}, by the It\^o formula,
\begin{align*}
&\psi(X^\e_t) - \int_0^t [{\ca{A}}_{r,\e}(X^\e_r)+{\ca{B}}_{r,\e}(X^\e_r)]dr\\
=&\psi(X^\e_t)-\int_0^t b^\e(r,X^\e_r)\cdot \nabla\psi(X^\e_r)dr-\frac12 \sum_{i,j=1}^d
\int_0^t a_{ij}^\e(r,X^\e_r) \partial_{ij}\psi(X^\e_r) dr\\
&- \int_0^t\int_E\int_{\R^d} \big[\psi(X^\e_{r}+h(s,z,x))-\psi(X^\e_{r})\big] F_s^\e(x,X^\e_r)
f_r(dx)\mu(dz)dr
\end{align*}
is a martingale, which implies the claim. We thus may write, for each $n\geq1$,
\begin{align*}
| \bb{E}[\ca{K}(X)] | \le& |\bb{E}[\ca{K}(X)]-\bb{E}[\tilde{\ca{K}}(X)] |
+  |\bb{E}[\tilde{\ca{K}}(X)]-\bb{E}[\tilde{\ca{K}}(X^{\e_n})] |\\
&+ |\bb{E}[\tilde{\ca{K}}(X^{\e_n})] - \bb{E}[\tilde{\ca{K}}_{\e_n}(X^{\e_n})] |  +
 |\bb{E}[\tilde{\ca{K}}_{\e_n}(X^{\e_n})] - \bb{E}[\ca{K}_{\e_n}(X^{\e_n})] |.
\end{align*}

We now study the four terms. We denote by $M$ a constant such that Supp $\psi\subset B(0,M)$.
We also define $\phi(z)=(2\pi)^{-d/2}e^{-|z|^2/2}$, so that $\phi_{\e}(z)=\e^{-d/2}\phi(\e^{-1/2}z)$.

{\it Step 1.} Here we prove that $\lim_{n\to\infty} \bb{E}[\tilde{\ca{K}}(X^{\e_n})]=\bb{E}[\tilde{\ca{K}}(X)]|$.
Since $X^{\e_n}$ goes in law to $X$ by construction, it suffices to verify that
$\tilde{\ca{K}}$ is bounded and a.s. continuous at $X$.

Since $\tilde{a}$, $\tilde{b}$  and $\tilde{h}$ are continuous in space and time, we easily deduce that
$(r,x)\mapsto \tilde {\ca A}_r\psi(x)$ and $(r,x)\mapsto \tilde {\ca{B}}_r\psi(x)$ are continuous and
bounded on $[0,T]\times \R^d$. For
$\tilde{\ca A}_r \psi(x)= \tilde b(r,x)\cdot \nabla\psi(x)
+\frac 1 2 \sum_{i,j} \tilde a_{ij}(r,x) \partial_{ij}\psi(x)$
this is obvious, and for $\tilde {\ca B}_r\psi(x)=\int_E [\psi(x+\tilde h(r,z,x))-\psi(x)]\mu(dz)
=\int_A[\psi(x+\tilde h(r,z,x))-\psi(x)]\mu(dz)$,
this follows from the Lebesgue theorem, because $\psi$ is bounded and $\mu(A)<\infty$.

We easily deduce that $\tilde{\ca{K}}$ is bounded, and that it is
continuous at each $\lambda\in\bb{D}([0,T],\R^d)$ which does not jump at
$s_1,\dots,s_k,s,t$. This is a.s. the case of $X$, see Lemma \ref{lem-tightness}.

{\it Step 2.} Here we check that
$\Delta_1:=|\bb{E}[\ca{K}(X)]-\bb{E}[\tilde{\ca{K}}(X)] | \leq C \rho$ for some constant $C$.
We have, since Supp $\psi \subset B(0,M)$,
\begin{align*}
|\ca{K}(\lambda)-\tilde{\ca{K}}(\lambda)|\leq& C \int_0^t [|\ca A_r\psi(\lambda_r)-\tilde{\ca A}_r\psi(\lambda_r)|
+|\ca B_r\psi(\lambda_r)-\tilde{\ca B}_r\psi(\lambda_r)|]dr \\
\leq & C \int_0^t \Big(|a(r,\lambda_r) - \tilde a(r,\lambda_r)|+ |b(r,\lambda_r) - \tilde b(r,\lambda_r)|\Big)
\bbd{1}_{\{|\lambda_r|<M\}} dr \\
&+ C \int_0^t \int_E |h(r,z,\lambda_r) - \tilde h(r,z,\lambda_r)| \mu(dz)dr.
\end{align*}
Using now that $\bbd{1}_{\{|x|<M\}} \leq C (1+|x|)^{-1}$ and that $\ca L(X_r)=f_r$ for each $r\in[0,T]$, we conclude
that
\begin{align*}
\Delta_1
\leq & C \int_0^t \int_{\R^d}\Big( \frac{|a(r,x) - \tilde a(r,x)|}{1+|x|}+
|b(r,x) - \tilde b(r,x)|\Big)f_r(dx)dr\\
&+ C \int_0^t \int_E \int_{\R^d} |h(r,z,x) - \tilde h(r,z,x)| f_r(dx)\mu(dz)dr.
\end{align*}
This is smaller than $C\rho$ by Lemma \ref{lem-continuous}.

{\it Step 3.} Now we verify that for all $n\geq 1$,
$\Delta_2^n=|\bb{E}[\tilde{\ca{K}}_{\e_n}(X^{\e_n})] - \bb{E}[\ca{K}_{\e_n}(X^{\e_n})]|\leq C\rho$.
As in Step 2,
\begin{align*}
\Delta_2^n
\le&  C \int_0^t \int_{\R^d} \Big(\frac{|a^{\e_n}(r, y)-\tilde{a}^{\e_n}(r,y)|}{1+|y|}+
|b^{\e_n}(r, y)-\tilde{b}^{\e_n}(r,y)|\Big) f_r^{\e_n}(y)dy \,dr\\
&+ C \, \int_0^t\int_E\int_{\R^d}\int_{\R^d}|h(r, z, x)-\tilde{h}(r,z,x)| \frac{\phi_{\e_n}(x-y)}{f^{\e_n}_t(y)}
\, f_r(dx) \,f_r^{\e_n}(y) dy \,\mu(dz) \,dr.
\end{align*}
Recalling (see Lemma \ref{regu-eq}) that $a^{\e_n}(r, y)f_r^{\e_n}(y)=\int_{\R^d}\phi_{\e_n}(x-y)a(r,x)f_r(dx)$,
that $\tilde a^{\e_n}(r, y)f_r^{\e_n}(y)=\int_{\R^d}\phi_{\e_n}(x-y)\tilde a(r,x)f_r(dx)$ and similar formulas
for $b^{\e_n}(r, y)f_r^{\e_n}(y)$ and $\tilde b^{\e_n}(r, y)f_r^{\e_n}(y)$, we find
\begin{align*}
\Delta_2^n
\le&  C \int_0^t \int_{\R^d}\int_{\R^d} \Big(\frac{|a(r, x)-\tilde{a}(r,x)|}{1+|y|}
+|b(r, x)-\tilde{b}(r,x)|\Big)\phi_{\e_n}(x-y) f_r(dx) dy \,dr\\
&+ C \, \int_0^t\int_E\int_{\R^d}\int_{\R^d}|h(r, z, x)-\tilde{h}(r,z,x)| \phi_{\e_n}(x-y)
\, f_r(dx) dy \,\mu(dz) \,dr.
\end{align*}
But $\int_{\R^d}\phi_{\e_n}(x-y)dy=1$ and,
since $\frac{1+|x|}{1+|y|}= 1+ \frac{|x|-|y|}{1+|y|}\leq 1+|x-y|
\leq 2+|x-y|^2$,
\begin{align*}
\int_{\R^d} \frac{(1+|x|)\phi_{\e_n}(x-y)dy}{1+|y|}\leq & \int_{\R^d} (2+ |x-y|^2)\phi_{\e_n}(x-y)dy
= 2 + d\e_n\leq 2+d.
\end{align*}
Consequently,
\begin{align*}
\Delta_2^n
\leq & C \int_0^t \int_{\R^d}\Big( \frac{|a(r,x) - \tilde a(r,x)|}{1+|x|}+
|b(r,x) - \tilde b(r,x)|\Big)f_r(dx)dr\\
&+ C \int_0^t \int_E \int_{\R^d} |h(r,z,x) - \tilde h(r,z,x)| f_r(dx)\mu(dz)dr,
\end{align*}
which is smaller than $C\rho$ by Lemma \ref{lem-continuous}.

{\it Step 4.} Finally, we check that $\lim_{n\to\infty}|\bb{E}[\tilde{\ca{K}}(X^{\e_n})] -
\bb{E}[\tilde{\ca{K}}_{\e_n}(X^{\e_n})]|=0$. We first observe that
$|\bb{E}[\tilde{\ca{K}}(X^{\e_n})] - \bb{E}[\tilde{\ca{K}}_{\e_n}(X^{\e_n})]| \le  C\, (I_n+J_n),$
where
\begin{align*}
I_n := \bb{E}\Big[\int_0^t |\tilde{\ca{A}}_{r,\e_n}\psi(X_r^{\e_n})-\tilde{\ca{A}}_{r}\psi(X_r^{\e_n})| dr
\Big]\quad\hbox{and}\quad
J_n := \bb{E}\Big[\int_0^t |\tilde{\ca{B}}_{r,\e_n}\psi(X_r^{\e_n})-\tilde{\ca{B}}_{r}\psi(X_r^{\e_n})|
dr  \Big].
\end{align*}
Since $\psi\in C^2_c(\R^d)$ and since $\ca{L}(X_r^{\e_n})=f_r^{\e_n}$, we have
\begin{align*}
I_n &\le C \int_0^t\int_{\R^d}\big(|\tilde{b}^{\e_n}(r,y)-\tilde{b}(r,y)|
+|\tilde{a}^{\e_n}(r,y)-\tilde{a}(r,y)|\big) \, f_r^{\e_n}(y)dy\, dr\\
&\leq C  \int_0^t\int_{\R^d}\int_{\R^d} \big(|\tilde{b}(r,x)-\tilde{b}(r,y)|
+|\tilde{a}(r,x)-\tilde{a}(r,y)|\big) \,\phi_{\e_n}(x-y) \, f_r(dx)\, dy dr.
\end{align*}
because $[\tilde b^{\e_n}(r, y)-\tilde b(r,y)]f_r^{\e_n}(y)
=\int_{\R^d}\phi_{\e_n}(x-y)\tilde b(r,x)f_r(dx) - \int_{\R^d} \phi_{\e_n}(x-y) \tilde b(r,y) f_r(dx)$,
with a similar formula concerning $\tilde a$. Using finally the
substitution $y=x+\sqrt{\e_n} u$, we find
$$I_n\leq C  \int_0^t\int_{\R^d}\int_{\R^d} \big(|\tilde{b}(r,x)-\tilde{b}(r,x+\sqrt{\e_n}u)|
+|\tilde{a}(r,x)-\tilde{a}(r,x+\sqrt{\e_n}u)|\big) \,\phi(u) \, f_r(dx)\, dy dr.
$$
Hence $\lim_n I_n=0$ by dominated convergence, since $\tilde a$ and $\tilde b$ are continuous and bounded.

By the same way, since
$f^{\e_n}_r(y)=\int_{\R^d}\phi_{\e_n}(x-y)f_r(dx)$,
\begin{align*}
J_n&=\bb{E}\Big[\int_0^t \Big|\int_E \int_{\R^d}\Big[\psi(X_r^{\e_n}+\tilde{h}(r,z,x))-\psi(X^{\e_n}_r)\Big]
\frac{\phi_{\e_n}(x-X_r^{\e_n})}{f_r^{\e_n}(X_r^{\e_n})} f_r(dx)\mu(dz) \\
&\hskip6cm - \int_E\Big[\psi(X_r^{\e_n}+\tilde{h}(r,z,X^{\e_n}_r))-\psi(X^{\e_n}_r)\Big]\mu(dz)\Big| dr \Big]\\
&=\bb{E}\Big[\int_0^t \Big|\int_E \int_{\R^d}\Big[\psi(X_r^{\e_n}+\tilde{h}(r,z,x))-
\psi(X_r^{\e_n}+\tilde{h}(r,z,X^{\e_n}_r))\Big]
\frac{\phi_{\e_n}(x-X_r^{\e_n})}{f_r^{\e_n}(X_r^{\e_n})} f_r(dx)\mu(dz)\Big|dr \Big]\\
&\leq C \bb{E}\Big[\int_0^t \int_E \int_{\R^d}\Big[ 1 \land \Big|\tilde{h}(r,z,x))-
\tilde{h}(r,z,X^{\e_n}_r)\Big|\Big]
\frac{\phi_{\e_n}(x-X_r^{\e_n})}{f_r^{\e_n}(X_r^{\e_n})} f_r(dx)\mu(dz)dr \Big]
\end{align*}
because $\psi$ and $\nabla \psi$ are bounded.
Using that $\ca L (X_r^{\e_n})=f_r^{\e_n}$, the substitution $y=x+\sqrt{\e_n} u$ and the fact that
$\tilde h(r,z,x)=0$ if $z\notin A$,
\begin{align*}
J_n&\leq C  \int_0^t \int_A \int_{\R^d}\int_{\R^d} \Big[1 \land |\tilde h(r,z,x)-\tilde h(r,z,y)|\Big]
\phi_{\e_n}(x-y) \, f_r(dx)\, dy \mu(dz) dr\\
&= C \int_0^t \int_A \int_{\R^d} \int_{\R^d} \Big[1 \land |\tilde h(r,z,x)-\tilde h(r,z,x+\sqrt{\e_n}u))|\Big]
\phi(u) \, f_r(dx)\, dy \mu(dz) dr.
\end{align*}
Hence $\lim_n J_n=0$ by dominated convergence, since $h$ is continuous in $x$ and since $\mu(A)<\infty$.

{\it Conclusion.} Gathering Steps 1, 2, 3 and 4, we find that $|\bb{E}[\ca{K}(X)]|\leq C \rho$.
Since $\rho$ can be chosen arbitrarily small, we conclude that $\bb{E}[\ca{K}(X)]=0$,
which completes the proof.
\epf

\section{Appendix}

\begin{proof}[Proof of Remark \ref{phiext}]
First, it is very easy, using only that $a$ and $b$ are locally bounded on $[0,T]\times \R^d$,
to show that $\ca A_t\varphi(x)$ is uniformly bounded as soon as $\varphi \in C^2_c(\R^d)$.
The case of $\ca B_t\varphi$ is more complicated. We consider
$\varphi \in C^2_c(\R^d)$ and $M>0$ such that Supp $\varphi\subset B(0,M)$ and we write
\begin{align*}
|\ca{B}_{t}\varphi (x)|\le & \bbd{1}_{\{|x|\le 2M\}} ||\nabla \varphi||_\infty \int_E|h(t,z,x)| \mu(dz)
+\bbd{1}_{\{|x|\ge 2M\}} \int_E|\varphi(x+h(t,z,x))|\mu(dz).
\end{align*}
We observe that $|\varphi(x+h(t,z,x))|\leq ||\varphi||_\infty\bbd{1}_{\{|x+h(t,z,x)|\le M\}}$ and that
\[
\bbd{1}_{\{|x|\ge 2M, |x+h(t,z,x)|\le M\}}\le \bbd{1}_{\{|x|\ge 2M, |h(t,z,x)|\ge |x|/2\}}
\le \bbd{1}_{\{|x|\ge 2M\}}\frac{2|h(t,z,x)|}{|x|}.
\]
Since $\int_E|h(t,z,x)|\mu(dz)\leq C(1+|x|)$ by assumption, we conclude that
$$
|\ca{B}_{t}\varphi (x)|\le \bbd{1}_{\{|x|\le 2M\}} C ||\nabla \varphi||_\infty(1+|x|)+ \bbd{1}_{\{|x|\ge 2M\}}
\frac{C  || \varphi||_\infty(1+|x|)}{|x|},
$$
which is bounded. We have proved point (i).

We next prove (ii). We put $\varphi(x)=(1+|x|^2)^{1/2}$, which satisfies
$$
\frac{1+|x|}2\leq \varphi(x)\leq 1+|x|,\quad |\nabla \varphi |\leq 1 \quad \hbox{and}
\quad|D^2 \varphi| \leq \frac C \varphi.
$$
We also introduce an increasing $C^2$ function $\chi : \R_+\mapsto \R_+$ such that
$\chi(r)=r$ for $r\in[0,1]$ and $\chi(r)=2$ for $r\geq 2$. We thus have
$$
r \land 1 \leq \chi(r)\leq 2(r\land 1), \quad
|\chi'(r)|\leq C \indiq_{\{r\leq 2\}} \quad \hbox{and}\quad |\chi''(r)|\leq C  \indiq_{\{1\leq r\leq 2\}}.
$$
We then set, for $n\geq 1$ and $x\in\R^d$, $\psi_n(x)=n \chi(\varphi(x)/n)$,
which satisfies
$$
\varphi \land n \leq \psi_n\leq 2(\varphi\land n), \quad |\nabla \psi_n|\leq C \indiq_{\{ \varphi \leq 2n\}}
\quad \hbox{and}\quad |D^2 \psi_n| \leq \frac C\varphi \indiq_{\{ \varphi \leq 2n\}}.
$$
Consequently, for all $s\in[0,T]$, since
$|b(s,\cdot)|\leq C \varphi$ and $|a(s,\cdot)|\leq C \varphi^2$ by Assumption \ref{Assum},
$$
|\ca A_s \psi_n |\leq |b(s,\cdot)||\nabla \psi_n| +  |a(s,\cdot)||D^2 \psi_n| \leq
C \varphi  \indiq_{\{ \varphi \leq 2n\}} \leq C [\varphi \land (2n)] \leq C \psi_n.
$$
We next claim that
\begin{equation}\label{claim}
\Delta_n(s,z,x)=|\psi_n(x+h(s,z,x))-\psi_n(x)|\leq C |h(s,z,x)| \frac{\psi_n(x)}{\varphi(x)}.
\end{equation}
First, if $\varphi(x)\leq 4n$, then we only use that $\nabla \psi_n$
is uniformly bounded to write $\Delta_n(s,z,x)\leq C |h(s,z,x)|$, whence the result because
$\psi_n(x)\geq \varphi(x)\land n \geq \varphi(x)/4$.
Second, if $\varphi(x)\geq 4n$ (whence $|x|\geq 4n-1\geq 3n$),
since $\psi_n$ is constant (with value $2n$) on $B(0,2n)^c$ and bounded on $\R^d$ by
$2n$, we can write
$\Delta_n(s,z,x)\leq 4n \indiq_{\{|x+h(s,z,x)|\leq 2n\}}\leq 4n \indiq_{\{|h(s,z,x)|\geq |x|/3\}} \leq 12n |h(s,z,x)| / |x|$.
But $12n =6\psi_n(x)$ and $|x| \geq \varphi(x)-1 \geq \varphi(x)/2$, whence the result.

We deduce from \eqref{claim}, using Assumption \ref{Assum}, that
$$
|\ca B_s \psi_n(x)| \leq C \frac{\psi_n(x)}{\varphi(x)} \int_E |h(s,z,x)|\mu(dz) \leq
C  \frac{\psi_n(x)}{\varphi(x)} (1+|x|)\leq
C \psi_n(x).
$$
Applying \eqref{IPDE} with the test function
$\psi_n - 2n \in C^2_c(\R^d)$, for which of course $(\ca A_s + \ca B_s)(\psi_n - 2n)=(\ca A_s + \ca B_s)\psi_n$,
and using that $f_0$ and $f_t$ are probability measures, we find
\begin{align*}
\int_{\R^d} \psi_n(x) f_t(dx)=&\int_{\R^d} \psi_n(x) f_0(dx)
+\int_0^t \int_{\R^d} (\ca A_s\psi_n(x)+ \ca B_s\psi_n(x)) f_s(dx) ds\\
\leq& \int_{\R^d}\psi_n(x)f_0(dx) + C \int_0^t\int_{\R^d} \psi_n(x)f_s(dx) ds.
\end{align*}
Since $f_0 \in \ca P_1(\R^d)$ by assumption and since $0\leq\psi_n(x)\leq 2|x|+2$,
$\sup_{n \geq 1}\int_{\R^d}\psi_n(x)f_0(dx) <\infty$.
We thus conclude, by the Gronwall Lemma, that
$\sup_{n \geq 1}\sup_{t \in [0,T]}\int_{\R^d} \psi_n(x)f_t(dx)<\infty$, which clearly implies that
$(f_t)_{t\in [0,T]} \in L^\infty([0,T],\ca P_1(\R^d))$, because $\lim_{n\to 0}\psi_n(x)=\varphi(x)\geq |x|$.

For point (iii),
we introduce a family of functions $\chi_n \in C^2_c(\R^d)$, for $n\geq 1$,
such that 
$\indiq_{\{|x|\leq n\}}\leq \chi_n(x) \leq \indiq_{\{|x|\leq n+1\}}$ 
and such that
$|D \chi_n (x)|+|D^2 \chi_n(x)|\leq C \bbd{1}_{\{|x| \in [n,n+1]\}}$.
We then consider $\varphi \in C^2(\R^d)$ as in the statement, i.e. such that
$(1+|x|)[|\varphi(x)|+|\nabla \varphi(x)|+|D^2\varphi(x)|]$ is bounded. Of course,
$\varphi\chi_n \in C^2_c(\R^d)$
for each $n\geq 1$, so that we can apply \eqref{IPDE}. We then let $n\to \infty$. Since $\varphi$ is bounded,
we obviously have $\lim_n \int_{\R^d} \varphi(x)\chi_n(x) f_t(dx) =\int_{\R^d} \varphi(x) f_t(dx)$.
Next, we want to prove that
$\lim_n \int_0^t \int_{\R^d} [\ca A_s (\varphi\chi_n)(x)+\ca B_s (\varphi\chi_n)(x)] f_s(dx)ds
= \int_0^t \int_{\R^d} [\ca A_s \varphi (x)+\ca B_s \varphi (x)] f_s(dx)ds$.
By dominated convergence and since $(f_t)_{t\in [0,T]} \in L^\infty([0,T],\ca P_1(\R^d))$
by (ii), it suffices to prove that for all $s\in[0,T], x\in \R^d$,

(a) $\sup_n|\ca A_s (\varphi\chi_n)(x)|\leq C(1+|x|)$,
(b) $\lim_n \ca A_s (\varphi\chi_n)(x)=\ca A_s \varphi (x)$,

(c) $\sup_n|\ca B_s (\varphi\chi_n)(x)|\leq C(1+|x|)$,
(d) $\lim_n \ca B_s (\varphi\chi_n)(x)=\ca B_s \varphi (x)$.

Point (a) is easy: since $|a(s,x)|+|b(s,x)|\leq C(1+|x|^2)$ by Assumption \ref{Assum}
and since $\chi_n,D\chi_n,D^2\chi_n$ are uniformly bounded,
$$
|\ca A_s (\varphi\chi_n)(x)|\leq C(1+|x|^2)(|D(\varphi\chi_n)(x)|+|D^2(\varphi\chi_n)(x)|)
\leq C(1+|x|^2)(|\varphi(x)|+|D\varphi(x)|+|D^2\varphi(x)|),
$$
which is bounded by $C(1+|x|)$ by assumption.
Point (b) is not hard, using that
$\lim_n \nabla (\varphi\chi_n)(x)=\nabla \varphi(x)$ and
$\lim_n \partial_{ij}(\varphi\chi_n)(x)=\partial_{ij}
\varphi(x)$ for each $x\in\R^d$.

Next, $\nabla (\varphi \chi_n)$ is uniformly bounded, so that
$|(\varphi\chi_n)(x+h(s,z,x))-(\varphi\chi_n)(x)| \leq C |h(s,z,x)|$ and thus
$|\ca B_s(\varphi\chi_n)(x)|\leq C \int_E |h(s,z,x)| \mu(dz)\leq C(1+|x|)$ by Assumption \ref{Assum},
whence (c).
Also, by dominated convergence, since  $\lim_n\chi_n(y)=1$ for all $y\in\R^d$,
$$
\lim_n \ca B_s(\varphi\chi_n)(x)=\lim_n\int_E [(\varphi\chi_n)(x+h(s,z,x))-(\varphi\chi_n)(x)]\mu(dz) =
\int_E [\varphi(x+h(s,z,x))-\varphi(x)]\mu(dz),
$$
which is nothing but $\ca B_s\varphi(x)$ as desired.
\end{proof}

{\bf Acknowledgments.}
We warmly thank Maxime Hauray for his help.

\bibliographystyle{siam} 
\bibliography{jump}
\end{document}